# Sampled-Data Observers for 1-D Parabolic PDEs with Non-Local Outputs

**Iasson Karafyllis[*], Tarek Ahmed-Ali[**] and Fouad Giri[**]**

[*]Dept. of Mathematics, National Technical University of Athens,
Zografou Campus, 15780, Athens, Greece,
email: iasonkar@central.ntua.gr; iasonkaraf@gmail.com

[**]Normandie UNIV, UNICAEN, ENSICAEN, LAC, 14000 Caen, France
email: tarek.ahmed-ali@ensicaen.fr; fouad.giri@unicaen.fr

**Abstract**

The present work provides a systematic approach for the design of sampled-data observers to a wide class of 1-D, parabolic PDEs with non-local outputs. The studied class of parabolic PDEs allows the presence of globally Lipschitz nonlinear and non-local terms in the PDE. Two different sampled-data observers are presented: one with an inter-sample predictor for the unavailable continuous measurement signal and one without an inter-sample predictor. Explicit conditions on the upper diameter of the (uncertain) sampling schedule for both designs are derived for exponential convergence of the observer error to zero in the absence of measurement noise and modeling errors. Moreover, explicit estimates of the convergence rate can be deduced based on the knowledge of the upper diameter of the sampling schedule. When measurement noise and/or modeling errors are present, Input-to-Output Stability (IOS) estimates of the observer error hold for both designs with respect to noise and modeling errors. The main results are illustrated by two examples which show how the proposed methodology can be extended to other cases (e.g., boundary point measurements).

## 1. Introduction

The problem of designing observers for PDEs has received a great deal of interest over the past decade, see e.g. [25,23,9,31,20,15,27,12]. Most existing works have been devoted to observer design for linear PDEs of parabolic and hyperbolic PDE type, using various design techniques including semigroup-based Luenberger method, modal decomposition, backstepping technique, Lyapunov stability, and LMIs. Observers for nonlinear PDEs have been proposed in e.g. [22,24,14,3]. Further interesting results have been reported on observer design for compound systems. In [22], backstepping observer design has been developed for linear ODE-PDE cascades with parabolic PDE. The result has been extended in [2] to cope with strict-feedback Lipschitz nonlinearities in the ODE part. The case of linear ODE-PDE cascades with first-order hyperbolic PDEs have been dealt with in [21]. Boundary observer design for linear PDE-ODEs (with hyperbolic PDE), has been considered in [32].

A common characteristic of all previously mentioned works is that it is supposed in all of them that the outputs are continuously accessible to measurements. However, usually only sampled (in time) measurements are available in practice. The observers that are based on continuous-time measurements can hardly achieve their theoretical performances if applied in the presence of measurement sampling. Therefore, an increasing research activity has recently been devoted to the problem of designing observers, both for ODEs and PDEs, that only require sampled (in time) output measurements. Most results on sampled-measurement observer design have been achieved for ODEs and only few for PDEs. In [11,30] sampled-output measurement observers of Luenberger-type have been designed for semilinear parabolic PDEs, using Lyapunov and LMIs tools. An extension of the results of [11] has been proposed in [4] to achieve larger sampling intervals by an inter-sample decreasing gain in the observer. Another extension has been proposed in [29] to solve the H1 filtering problem in the case of a reaction-diffusion system with sampled measurements. Sampled-output measurement observer design for ODE-PDE cascades has also recently been investigated in [5,6] for a class of parabolic PDEs.

Observers for PDEs (with sampled-measurements or not) can be divided in two classes depending on the definition of the system output (i.e. the signal that provides the information used in the observer). The first class includes observers which require as output a vector formed by a finite number of measurements provided physical sensors, each one providing the measurements of the PDE state at the point where it is placed on the domain. The observers in this class are observers which use point, local or sampled in space measurements; see for instance [3,7,11]. A special case is the case where the sensors are placed at the domain boundary (providing the PDE state or its time-derivative at the boundary); the resulting observer is termed as boundary observer. This is the case of most existing observers e.g. [9,31,15,20,14,21,4]. The second class of observers is constituted of those observers that deal with non-local outputs. The outputs in this case are functionals of the state defined so that the output values are weighted averages of the state values over the domain (see [8,10,13]).

A crucial characteristic of sampled-data observers (those for ODEs and PDEs) is the way the innovation term depends on the (sampled) output measurements. Typically, the existing observers are classified in two categories. The first category is one where the observer involves a Zero-Order-Hold (ZOH) innovation term. Accordingly, the innovation term is updated at each sampling time, using the output measurement sample, and kept unchanged between two successive sampling times. Observers for PDEs belonging to this category are those proposed in [7,8,11,10,10,13]. Many more observers of this type have been proposed for ODEs. The second category of sampled-data observers is characterized by innovation terms that are continuously updated by using output prediction between two successive sampling times. It turns out that the observers of this type include an additional crucial component referred to as inter-sample output predictor. The approach was first proposed in [17] and several sampled-data observers involving inter-sample output predictors have been developed for different classes of ODEs, e.g. [1,16]. The observers with inter-sample predictors are usually more convenient in practice, since they allow less frequent sampling compared to ZOH-innovation observers.

The problem of designing sampled-data observers with inter-sample predictors for PDEs has yet to be solved. Indeed, the extension of the approach in [17] to parabolic PDEs is challenging, since the time derivative of a functional may not be expressed in a way that allows the derivation of bounds for the time derivative. More specifically, the derivation of bounds for the time derivative of a functional may involve quantities which cannot be estimated (e.g., the second spatial derivative of the observer error).

The present paper shows that the approach in [17] can be applied to a wide class of 1-D, parabolic PDEs with non-local outputs. Furthermore, the studied class of parabolic PDEs allows the presence of globally Lipschitz nonlinearities and non-local terms in the PDE. Moreover, the sampling schedule is not required to be exactly known. Applying a small-gain approach and modifying the construction of the inter-sample predictor, we are in a position to derive explicit conditions on the upper diameter of the uncertain sampling schedule which guarantee exponential convergence of the observer error to zero in the absence of noise and modeling errors (Theorem 2.2). The results are compared with the ZOH-innovation observer design (Theorem 2.3). Explicit estimates of the convergence rate can be deduced in both cases. On the other hand, when noise and/or modeling errors are present, we are in a position to guarantee Input-to-Output Stability (IOS) estimates of the observer error with respect to noise and modeling errors. Two examples are presented in Section 3 of the present paper, which aim to show how easily we can design sampled-data observers. Interestingly, the examples also illustrate additional facts:

- the first example (Example 3.1) shows that there exist 1-D, parabolic PDEs which allow the diameter of the sampling schedule to be arbitrarily large for the sampled-data observer with the inter-sample predictor, while the sampled-data observer without the inter-sample predictor requires a sufficiently frequent sampling schedule, and
- the second example (Example 3.2) shows that our proposed methodology can be used also when the measurements are local outputs (e.g., boundary point measurements) and can guarantee stronger estimates of the observer error (e.g., estimates in the spatial sup norm).

The structure of the paper is as follows. In Section 2, we describe the class of systems which are studied and the construction of the sampled-data observers. The main results (Theorem 2.2 and Theorem 2.3) are also stated in Section 2. Section 3 contains the two examples that illustrate the use of the obtained results for the design of sampled-data observers. The proofs of the main results are contained in Section 4. Finally, the concluding remarks of the present work are given in Section 5.

**Notation:** Throughout the paper, we adopt the following notation.

* $\Re_+ := [0, +\infty)$. Let $u : \Re_+ \times [0,1] \to \Re$ be given. We use the notation $u[t]$ to denote the profile at certain $t \geq 0$, i.e., $(u[t])(x) = u(t,x)$ for all $x \in [0,1]$. $L^2(0,1)$ denotes the equivalence class of measurable functions $f : [0,1] \to \Re$ for which $\|f\| = \left( \int_0^1 |f(x)|^2 dx \right)^{1/2} < +\infty$. For an integer $k \geq 1$, $H^k(0,1)$ denotes the Sobolev space of functions in $L^2(0,1)$ with all its weak derivatives up to order $k \geq 1$ in $L^2(0,1)$.

* Let $S \subseteq \Re^n$ be an open set and let $A \subseteq \Re^n$ be a set that satisfies $S \subseteq A \subseteq cl(S)$. By $C^0(A; \Omega)$, we denote the class of continuous functions on $A$, which take values in $\Omega \subseteq \Re^m$. By $C^k(A; \Omega)$, where $k \geq 1$ is an integer, we denote the class of functions on $A \subseteq \Re^n$, which takes values in $\Omega \subseteq \Re^m$ and has continuous derivatives of order $k$. When $\Omega = \Re$ then we write $C^0(A)$ or $C^k(A)$. For $f \in C^0([0,1])$ the sup norm is defined by $\|f\|_\infty = \sup_{0 \leq z \leq 1}(|f(z)|) < +\infty$.

* A continuous function $f : [0,1] \to \Re$ is called piecewise $C^1$ on $[0,1]$ and we write $f \in PC^1([0,1])$, if the following properties hold: (i) for every $x \in [0,1)$ the limit $\lim_{h \to 0^+} \left( h^{-1}(f(x+h) - f(x)) \right)$ exists and is finite, (ii) for every $x \in (0,1]$ the limit $\lim_{h \to 0^-} \left( h^{-1}(f(x+h) - f(x)) \right)$ exists and is finite, (iii) there exists a set $J \subset (0,1)$ of finite cardinality, where $\frac{df}{dx}(x) = \lim_{h \to 0^-} \left( h^{-1}(f(x+h) - f(x)) \right) = \lim_{h \to 0^+} \left( h^{-1}(f(x+h) - f(x)) \right)$ holds for $x \in (0,1) \setminus J$, and (iv) the mapping $(0,1) \setminus J \ni x \to \frac{df}{dx}(x) \in \Re$ is continuous.

* For a vector $x \in \Re^n$ we denote by $|x|$ its usual Euclidean norm, by $x'$ its transpose. For a real matrix $A \in \Re^{n \times m}$, $A' \in \Re^{m \times n}$ denotes its transpose and $|A| := \sup\{|Ax| ; x \in \Re^n, |x| = 1\}$ is its induced matrix norm. For two positive integers $i, j > 0$, $\delta_{i,j}$ denotes the Kronecker delta, i.e., $\delta_{i,j} = 1$ if $i = j$ and $\delta_{i,j} = 0$ if $i \neq j$.

## 2. System Description and Main Results

*2.A. System Description*

Consider the Sturm-Liouville (SL) operator $B : D \to L^2(0,1)$ with

$$(Bu)(x) = -p \frac{d^2 u}{dx^2}(x) + q(x)u(x), \text{ for } x \in (0,1) \text{ and } u \in D \tag{2.1}$$

where $p > 0$ is a constant, $q \in C^0([0,1])$,

$$D := \left\{ u \in H^2(0,1) : a_0 u(0) + b_0 \frac{du}{dx}(0) = a_1 u(1) + b_1 \frac{du}{dx}(1) = 0 \right\}, \quad (2.2)$$

with $a_0, b_0, a_1, b_1 \in \Re$ being constants with $a_0^2 + b_0^2 > 0$, $a_1^2 + b_1^2 > 0$. Let $\lambda_1 < \lambda_2 < ... < \lambda_n < ...$ and $\phi_n \in C^2([0,1]) \cap D$ ($n = 1, 2, ...$) (with $\|\phi_n\| = 1$) be the eigenvalues and the eigenfunctions of the SL operator $B : D \to L^2(0,1)$. In this work we make the following assumption.

**(H1):** *The SL operator $B : D \to L^2(0,1)$ defined by (2.1), (2.2), where $a_0, b_0, a_1, b_1 \in \Re$ are constants with $a_0^2 + b_0^2 > 0$, $a_1^2 + b_1^2 > 0$, satisfies*

$$\sum_{n=M}^{\infty} \lambda_n^{-1} \max_{0 \leq x \leq 1} (|\phi_n(x)|) < +\infty, \text{ for certain } M > 0 \text{ with } \lambda_M > 0 \quad (2.3)$$

It is important to notice that the validity of Assumption (H1) can be verified without the knowledge of the eigenvalues and eigenfunctions of the SL operator $B$. In effect, it is shown in [26] that Assumption (H1) holds automatically, provided that $b_0, a_1, b_1 \geq 0$, $a_0 \leq 0$.

With the help of the SL operator $B : D \to L^2(0,1)$ defined by (2.1), (2.2) under Assumption (H1), we are in a position to describe the system that we study. We consider the observer design problem for the system that is described by the parabolic PDE:

$$\frac{\partial u}{\partial t}(t,x) = p \frac{\partial^2 u}{\partial x^2}(t,x) - q(x) u(t,x) + K(u[t])(x) + g(x, \bar{P}u[t]) + v(t,x), \text{ for } (t,x) \in (0,+\infty) \times (0,1) \quad (2.4)$$

$$a_0 u(t,0) + b_0 \frac{\partial u}{\partial x}(t,0) = a_1 u(t,1) + b_1 \frac{\partial u}{\partial x}(t,1) = 0, \text{ for } t \geq 0 \quad (2.5)$$

where the nature of all mappings appearing in (2.4) are explained by the following assumptions.

**(H2)** *The following regularity requirements hold.*
- $g \in C^1([0,1] \times \Re^k)$, $K : C^0([0,1]) \to C^0([0,1])$ *are continuous mappings with $K(u) \in D$ for all $u \in C^0([0,1])$, for which there exists a constant $\bar{L} > 0$ such that the inequalities $\|K(u) - K(w)\|_\infty \leq \bar{L} \|u - w\|_\infty$, $\max_{0 \leq x \leq 1}(|g(x,\beta) - g(x,\gamma)|) \leq \bar{L}|\beta - \gamma|$, $\|BK(u)\| \leq \bar{L} \|u\|_\infty$ hold for all $u, w \in C^0([0,1])$, $\beta, \gamma \in \Re^k$,*
- $\bar{P} : C^0([0,1]) \to \Re^k$ *is a compatible operator with $B$, i.e., a continuous linear operator for which there exists a continuous operator $S : C^0([0,1]) \to \Re^k$ such that $\bar{P} Bu = -Su$ for all $u \in C^2([0,1]) \cap D$ (see [19]),*
- *the mapping $f : L^2(0,1) \to L^2(0,1)$ defined by*

$$f(u)(x) = K(u)(x) + g(x, \bar{P}u), \text{ for } x \in [0,1] \text{ and } u \in L^2(0,1) \quad (2.6)$$

*is a continuous mapping for which there exists a constant $R \geq 0$ such that the following global Lipschitz inequality holds:*

$$\|f(u) - f(w)\| \leq R \|u - w\|, \text{ for all } u, w \in L^2(0,1) \quad (2.7)$$

**(H3)** *There exist mappings $v_1, v_2 \in C^0(\Re_+ \times [0,1])$ for which the mapping $(0,+\infty) \times [0,1] \ni (t,x) \to \frac{\partial v_1}{\partial t}(t,x)$ is continuous with $v_1[t] \in PC^1([0,1])$ for all $t \geq 0$, $v_2[t] \in D$ and $\sup_{\tau \in (0,t)}(\|Bv_2[\tau]\|) < +\infty$ for all $t > 0$ and such that the (distributed) input $v : \Re_+ \times [0,1] \to \Re$ satisfies*

$$v[t] = v_1[t] + v_2[t], \text{ for } t \geq 0, \ x \in [0,1] \quad (2.8)$$

Using Corollary 4.6 in [19], we are in a position to guarantee that when Assumptions (H1), (H2), (H3) hold then for every $u_0 \in D$ there exists a unique mapping $u \in C^0(\Re_+ \times [0,1]) \cap C^1((0,+\infty) \times [0,1])$ satisfying $u[t] \in C^2([0,1])$ for all $t > 0$, (2.4), (2.5) and the initial condition

$$u[0] = u_0 \tag{2.9}$$

The measured outputs are non-local and are described by the equations

$$y_i(t) = \xi_i(t_j) + \int_0^1 k_i(x) u(t_j, x) dx, \text{ for } t \in [t_j, t_{j+1}), \ i = 1,...,m \text{ and } j = 0,1,2,... \tag{2.10}$$

where $k_i \in L^2(0,1)$ ($i = 1,...,m$) are the output kernels, $\xi = (\xi_1,...,\xi_m)': \Re_+ \to \Re^m$ is a bounded mapping (the measurement error) and $\{t_j \geq 0 : j = 0,1,2,...\}$ (the sampling schedule) is an increasing sequence of times (the sampling times) with $t_0 = 0$ and $\lim_{j \to +\infty}(t_j) = +\infty$. It is clear from (2.10) that the output $y(t) = (y_1(t),..., y_m(t))' \in \Re^m$ is finite-dimensional and sampled. Moreover, we won't assume that the sampling schedule is known (uncertain sampling schedule).

**Remark 2.1:** **(a)** The parabolic PDE system (2.3), (2.4) is a special case where the SL operator has a constant diffusion coefficient and zero advection terms. However, the general case where $B: D \to L^2(0,1)$ is defined by (2.2) and

$$(Bu)(x) = -\frac{1}{r(x)} \frac{\partial}{\partial x}\left( p(x) \frac{\partial u}{\partial x}(t,x) \right) + \frac{q(x)}{r(x)} u(t,x), \text{ for } x \in (0,1) \text{ and } u \in D \tag{2.11}$$

with $p, r \in C^2([0,1]; (0,+\infty))$, can be reduced to the case (2.1) by means of the transformation:

$$\xi = \sqrt{\varepsilon} \int_0^x \sqrt{\frac{p(s)}{r(s)}} ds \ , \text{ where } \varepsilon = \left( \int_0^1 \sqrt{\frac{r(s)}{p(s)}} ds \right)^{-2} \text{ and } U(t,\xi) = (r(x) p(x))^{1/4} u(t,x).$$

**(b)** The PDE system (2.3), (2.4) is allowed to contain nonlinear and non-local terms in the PDE.

*2.B. The Sampled-Data Observer with Inter-Sample Predictor*

The first proposed sampled-data observer consists of two components: the continuous-time observer and the inter-sample predictor. We start with the continuous-time observer.

Let $N \geq 1$ be an integer with $\lambda_{N+1} > 0$, $c_i \in D$ ($i = 1,...,m$) be a given set of functions and $L = \{L_{i,j} : i = 1,...,N, j = 1,...,m\} \in \Re^{N \times m}$ be a given matrix. Define the real matrix $A \in \Re^{N \times N}$ by

$$A_{i,j} = -\lambda_i \delta_{i,j} + L_{i,1} c_{1,j} + ... + L_{i,m} c_{m,j} \text{ for } i, j = 1,..., N \tag{2.12}$$

where

$$c_{i,j} = \int_0^1 c_i(s) \phi_j(s) ds, \text{ for } i = 1,...,m, \ j = 1,2,... \tag{2.13}$$

Also define

$$l_i(x) = \sum_{n=1}^N \phi_n(x) L_{n,i} \text{ for } i = 1,...,m \text{ and } x \in [0,1]. \tag{2.14}$$

The continuous-time observer is described by the following equations:

$$\frac{\partial w}{\partial t}(t,x) = p \frac{\partial^2 w}{\partial x^2}(t,x) - q(x) w(t,x) + K(w[t])(x) + g(x, \bar{P}w[t]) + \tilde{v}(t,x) + \sum_{i=1}^m l_i(x) \left( \int_0^1 c_i(s) w(t,s) ds - \zeta_i(t) \right)$$

$$\text{for all } (t,x) \in (0,+\infty) \times (0,1) \tag{2.15}$$

$$a_0 w(t,0) + b_0 \frac{\partial w}{\partial x}(t,0) = a_1 w(t,1) + b_1 \frac{\partial w}{\partial x}(t,1) = 0, \text{ for } t \geq 0 \tag{2.16}$$

where the distributed observer state $w[t]$ is to be used to approximate the state $u[t]$ and the additional observer states $\zeta_i(t)$ ($i=1,...,m$) are to be used to approximate the (unavailable) continuous signals $\int_0^1 c_i(x)u(t,x)dx$ ($i=1,...,m$). The distributed input $\tilde{v} \in C^0(\Re_+ \times [0,1])$ is assumed to satisfy Assumption (H3) and ideally it would be equal to $v$. However, we do not assume that $\tilde{v}$ coincides with $v$, in order to allow the expression of the effect of possible modeling errors.

The evolution of the observer states $\zeta_i(t)$ ($i=1,...,m$) is determined by the inter-sample predictor, which is described next. Ideally, we would like to have an inter-sample predictor for the output signals. However (by virtue of (2.4), (2.6)) the nominal output signals (i.e. without measurement noise) satisfy the following differential equations for $t>0$, $i=1,...,m$:

$$\frac{d}{dt}\int_0^1 k_i(x)u(t,x)dx = p\int_0^1 k_i(x)\frac{\partial^2 u}{\partial x^2}(t,x)dx - \int_0^1 q(x)k_i(x)u(t,x)dx + \int_0^1 k_i(x)(f(u[t]))(x)dx + \int_0^1 k_i(x)v(t,x)dx$$

and it becomes clear that the right-hand side of the above differential equation cannot be bounded by an estimate that involves the state norm $\|u[t]\|$. Therefore, we select $c_i \in D$ ($i=1,...,m$) to approximate closely the output kernels $k_i \in L^2(0,1)$ ($i=1,...,m$) and instead of using an inter-sample predictor for the output signals $\int_0^1 k_i(x)u(t,x)dx$ ($i=1,...,m$), we use the inter-sample predictor for the signals $\int_0^1 c_i(x)u(t,x)dx$ ($i=1,...,m$). Notice that the (unavailable) continuous signals $\int_0^1 c_i(x)u(t,x)dx$ ($i=1,...,m$) satisfy the following differential equations for $t>0$, $i=1,...,m$:

$$\frac{d}{dt}\int_0^1 c_i(x)u(t,x)dx = \int_0^1 \left(p\frac{d^2c_i}{dx^2}(x) - q(x)c_i(x)\right)u(t,x)dx + \int_0^1 c_i(x)(f(u[t]))(x)dx + \int_0^1 c_i(x)v(t,x)dx \quad (2.17)$$

Indeed, integrating by parts we get $\int_0^1 c_i(x)\frac{\partial^2 u}{\partial x^2}(t,x)dx = \psi_i(t) + \int_0^1 \frac{d^2c_i}{dx^2}(x)u(t,x)dx$, where

$\psi_i(t) = c_i(1)\frac{\partial u}{\partial x}(t,1) - c_i(0)\frac{\partial u}{\partial x}(t,0) - \frac{dc_i}{dx}(1)u(t,1) + \frac{dc_i}{dx}(0)u(t,0)$ for all $t>0$, $i=1,...,m$. The facts that $c_i \in D$ ($i=1,...,m$) and $a_0^2 + b_0^2 > 0$, $a_1^2 + b_1^2 > 0$, in conjunction with (2.2) and (2.5) implies that $\psi_i(t) = 0$ for all $t>0$, $i=1,...,m$. Equalities (2.17) follow from (2.4), (2.6) and the fact that $\int_0^1 c_i(x)\frac{\partial^2 u}{\partial x^2}(t,x)dx = \int_0^1 \frac{d^2c_i}{dx^2}(x)u(t,x)dx$ for $t>0$, $i=1,...,m$.

The inter-sample predictor replaces the (unavailable) state $u$ in (2.17) by its estimate $w$ and tries to approximate the unavailable signals $\int_0^1 c_i(x)u(t,x)dx$ ($i=1,...,m$). We get:

$$\zeta_i(t_j) = y(t_j) - \int_0^1 (k_i(x) - c_i(x))w(t_j,x)dx, \text{ for } i=1,...,m, \; j=0,1,2,... \quad (2.18)$$

$$\dot{\zeta}_i(t) = \int_0^1 \left(p\frac{d^2c_i}{dx^2}(x) - q(x)c_i(x)\right)w(t,x)dx + \int_0^1 c_i(x)(f(w[t]))(x)dx + \int_0^1 c_i(x)\tilde{v}(t,x)dx,$$
$$\text{for } t \in [t_j, t_{j+1}), \; i=1,...,m \text{ and } j=0,1,2,... \quad (2.19)$$

It should be noticed that Assumptions (H1), (H2), (H3) guarantee that for every $u_0, w_0 \in D$, for every input $\tilde{v} \in C^0(\Re_+ \times [0,1])$ that satisfies Assumption (H3) and for every increasing sequence of times $\{t_j \geq 0 : j = 0,1,2,...\}$ with $t_0 = 0$ and $\lim_{j \to +\infty}(t_j) = +\infty$, there exist unique mappings $w \in C^0(\Re_+ \times [0,1]) \cap C^1(\bar{I} \times [0,1])$ and $\zeta_i : \Re_+ \to \Re$ being right-continuous with $\zeta_i \in C^1(\bar{I})$ ($i = 1,...,m$), where $\bar{I} := \Re_+ \setminus \{t_j : j = 0,1,2,...\}$, satisfying $w[t] \in C^2([0,1])$ for all $t > 0$, $w[0] = w_0$, (2.18), (2.19) for all $t \in \bar{I}$, (2.15) for all $(t,x) \in \bar{I} \times (0,1)$ and (2.16). Indeed, the solution may be constructed by a step-by-step procedure. To see this, take any $j = 0,1,2,...$ and suppose that the solution $w$ on $t \in [0,t_j]$ and $\zeta_i$ ($i = 1,...,m$) on $t \in [0,t_j)$ (when $j > 0$) is already known. First use (2.18) in order to get the values of $\zeta_i$ ($i = 1,...,m$) for $t = t_j$. Notice that by virtue of Assumptions (H1), (H2), (H3) and the fact that $l_i \in D$ for $i = 1,...,m$ (recall (2.14)), all assumptions of Theorem 4.5 in [19] hold for system (2.15), (2.16), (2.19) on the interval $t \in [t_j, t_{j+1}]$. Therefore, using Theorem 4.5 in [19], we obtain the solution $w$ on $t \in [t_j, t_{j+1}]$ and $\zeta_i$ ($i = 1,...,m$) on $t \in [t_j, t_{j+1})$.

*2.C. The Sampled-Data Observer without the Inter-Sample Predictor*

The second proposed sampled-data observer is simpler than the first sampled-data observer since the inter-sample predictor is not used. The observer is described by (2.16) and the following equation:

$$\frac{\partial w}{\partial t}(t,x) = p\frac{\partial^2 w}{\partial x^2}(t,x) - q(x)w(t,x) + K(w[t])(x) + g(x, \bar{P}w[t]) + \tilde{v}(t,x) + \sum_{i=1}^{m} l_i(x)\left(\int_0^1 k_i(s)w(t_j,s)ds - y_i(t_j)\right)$$

for $(t,x) \in (t_j, t_{j+1}) \times (0,1)$ and $j = 0,1,2,...$ (2.20)

where (again) the distributed observer state $w[t]$ is to be used to approximate the state $u[t]$. Similarly to the first observer, the distributed input $\tilde{v} \in C^0(\Re_+ \times [0,1])$ is assumed to satisfy Assumption (H3) and ideally it would be equal to $v$. However, we do not assume that $\tilde{v}$ coincides with $v$, in order to allow the expression of the effect of possible modeling errors.

It should be noticed that Assumptions (H1), (H2), (H3) guarantee that for every $u_0, w_0 \in D$, for every input $\tilde{v} \in C^0(\Re_+ \times [0,1])$ that satisfies Assumption (H3) and for every increasing sequence of times $\{t_j \geq 0 : j = 0,1,2,...\}$ with $t_0 = 0$ and $\lim_{j \to +\infty}(t_j) = +\infty$, there exists a unique mapping $w \in C^0(\Re_+ \times [0,1]) \cap C^1(\bar{I} \times [0,1])$, where $\bar{I} := \Re_+ \setminus \{t_j : j = 0,1,2,...\}$, satisfying $w[t] \in C^2([0,1])$ for all $t > 0$, $w[0] = w_0$, (2.20) for all $(t,x) \in \bar{I} \times (0,1)$ and (2.16). Indeed, the solution may be constructed by a step-by-step procedure. To see this, take any $j = 0,1,2,...$ and suppose that the solution $w$ on $t \in [0,t_j]$ is already known. Notice that by virtue of Assumptions (H1), (H2), (H3) and the fact that $l_i \in D$ for $i = 1,...,m$ (recall (2.14)), all assumptions of Theorem 4.5 in [19] hold for system (2.20), (2.16) on the interval $t \in [t_j, t_{j+1}]$. Therefore, using Theorem 4.5 in [19], we obtain the solution $w$ on $t \in [t_j, t_{j+1}]$.

*2.D. Main Results*

We are now in a position to give conditions on the nonlinear term $f$ and the sampling schedule $\{t_j \geq 0 : j = 0,1,2,...\}$ that guarantee convergence of the estimation error $e[t] = w[t] - u[t]$ in the $L^2$ norm. Indeed, it is shown that the proposed sampled-data observers work provided that (i) the strength of the nonlinear term (i.e., the constant $R \geq 0$ appearing in (2.7)) is sufficiently small, (ii) the functions $c_i \in D$ ($i = 1,...,m$) approximate closely the output kernels $k_i \in L^2(0,1)$ ($i = 1,...,m$), and (iii) the sampling schedule is sufficiently frequent. It should be noticed that since $D$ (defined by (2.2)) is dense in $L^2(0,1)$, the close approximation of the output kernels $k_i \in L^2(0,1)$ ($i = 1,...,m$) by functions $c_i \in D$ ($i = 1,...,m$) is not a problem.

**Theorem 2.2 (Sampled-Data Observer Design with Inter-Sample Predictor in presence of measurement errors):** *Consider system (2.4), (2.5) under Assumptions (H1), (H2) with output given by (2.10), where $k_i \in L^2(0,1)$ ($i=1,...,m$) are the output kernels, $\xi: \Re_+ \to \Re^m$ is the measurement error and $\{t_j \geq 0: j=0,1,2,...\}$ is the sampling schedule. Let $N \geq 1$ be an integer with $\lambda_{N+1} > 0$, let $c_i \in D$ ($i=1,...,m$) be given functions and let $L = \{L_{i,j}: i=1,...,N, j=1,...,m\} \in \Re^{N \times m}$ be a matrix so that $A \in \Re^{N \times N}$ defined by (2.12) is a Hurwitz matrix. Define $l_i \in D$ for $i=1,...,m$ by means of (2.14) and*

$$K := \left( \sum_{j=N+1}^{\infty} c_{1,j}^2 + ... + \sum_{j=N+1}^{\infty} c_{m,j}^2 \right)^{1/2} \tag{2.21}$$

*Let $P \in \Re^{N \times N}$ be a symmetric, positive definite matrix with $P \geq I$ for which there exists a constant $\sigma > 0$ such that $PA + A'P \leq -2\sigma P$. Suppose that there exist constants $Q \geq 2$ with $Q > 2|L'PL|K^2/(\sigma \lambda_{N+1})$, $h > 0$, $\kappa \in [0, \mu)$, where $\mu := (H(Q) + 2\lambda_{N+1})/4$ and $H(Q) := 2\sigma - \lambda_{N+1} - \sqrt{(2\sigma - \lambda_{N+1})^2 + 16Q^{-1}|L'PL|K^2}$, such that the following small-gain condition holds:*

$$\Omega := \gamma \left( R + \exp(\kappa h) \sum_{i=1}^{m} \|l_i\| \left( \left( \left\| p \frac{d^2 c_i}{dx^2} - q c_i \right\| + R \|c_i\| \right) h + \|k_i - c_i\| \right) \right) < 1 \tag{2.22}$$

*where $\gamma := \sqrt{\dfrac{\tilde{g}}{2(\mu - \kappa)}}$ and $\tilde{g} := \max\left( \dfrac{4|P|}{4\sigma + H(Q)}, \dfrac{Q}{2\lambda_{N+1}} \right)$. Then for every $u_0, w_0 \in D$, for every bounded $\xi: \Re_+ \to \Re^m$, for every inputs $v, \tilde{v} \in C^0(\Re_+ \times [0,1])$ that satisfy Assumption (H3) and for every increasing sequence of times $\{t_j \geq 0: j=0,1,2,...\}$ with $t_0 = 0$, $\lim_{j \to +\infty} (t_j) = +\infty$ and $\sup_{j \geq 0}(t_{j+1} - t_j) \leq h$, the unique solution of the problem (2.4), (2.5), (2.10), (2.15), (2.16), (2.18), (2.19) with initial condition given by (2.9) and $w[0] = w_0$ satisfies the following estimate for $t \geq 0$ for the observer error $e[t] = w[t] - u[t]$:*

$$\|e[t]\| \leq (1-\Omega)^{-1} \sqrt{\max\left(|P|, \frac{Q}{2}\right)} \exp(-\kappa t) \|e[0]\| + \exp(\kappa h)(1-\Omega)^{-1} \gamma \sum_{i=1}^{m} \|l_i\| \sup_{0 \leq s \leq t} \left( |\xi_i(s)| \exp(-\kappa(t-s)) \right)$$

$$+ (1-\Omega)^{-1} \gamma \left( 1 + h \exp(\kappa h) \sum_{i=1}^{m} \|l_i\| \|c_i\| \right) \sup_{0 \leq s \leq t} \left( \|v[s] - \tilde{v}[s]\| \exp(-\kappa(t-s)) \right) \tag{2.23}$$

**Theorem 2.3 (Sampled-Data Observer Design without Inter-Sample Predictor in presence of measurement errors):** *Consider system (2.4), (2.5) under Assumptions (H1), (H2) with output given by (2.10), where $k_i \in L^2(0,1)$ ($i=1,...,m$) are the output kernels, $\xi: \Re_+ \to \Re^m$ is the measurement error and $\{t_j \geq 0: j=0,1,2,...\}$ is the sampling schedule. Let $N \geq 1$ be an integer with $\lambda_{N+1} > 0$, let $c_i \in D$ ($i=1,...,m$) be given functions and let $L = \{L_{i,j}: i=1,...,N, j=1,...,m\} \in \Re^{N \times m}$ be a matrix so that $A \in \Re^{N \times N}$ defined by (2.12) is a Hurwitz matrix. Define $l_i \in D$ for $i=1,...,m$ by means of (2.14) and $K \geq 0$ by (2.21). Let $P \in \Re^{N \times N}$ be a symmetric, positive definite matrix with $P \geq I$ for which there exists a constant $\sigma > 0$ such that $PA + A'P \leq -2\sigma P$. Suppose that there exist constants $Q \geq 2$ with $Q > 2|L'PL|K^2/(\sigma \lambda_{N+1})$, $h > 0$, $\kappa \in [0,\mu)$, where $H(Q) := 2\sigma - \lambda_{N+1} - \sqrt{(2\sigma - \lambda_{N+1})^2 + 16Q^{-1}|L'PL|K^2}$ and $\mu := (H(Q) + 2\lambda_{N+1})/4$, such that the following small-gain condition holds:*

$$\Omega := \gamma \left( R + \exp(\kappa h) \sum_{i=1}^{m} \|l_i\| \left( \left( \left\| p \frac{d^2 c_i}{dx^2} - q c_i \right\| + R \|c_i\| + \sum_{r=1}^{m} \left| \int_0^1 c_i(x) l_r(x) dx \right| \|k_r\| \right) h + \|k_i - c_i\| \right) \right) < 1 \tag{2.24}$$

where $\gamma := \sqrt{\dfrac{\tilde{g}}{2(\mu - \kappa)}}$ and $\tilde{g} := \max\left(\dfrac{4|P|}{4\sigma + H(Q)}, \dfrac{Q}{2\lambda_{N+1}}\right)$. *Then for every $u_0, w_0 \in D$, for every bounded $\xi: \Re_+ \to \Re^m$, for every inputs $v, \tilde{v} \in C^0(\Re_+ \times [0,1])$ that satisfy Assumption (H3) and for every increasing sequence of times $\{t_j \geq 0 : j = 0,1,2,\ldots\}$ with $t_0 = 0$, $\lim_{j \to +\infty}(t_j) = +\infty$ and $\sup_{j \geq 0}(t_{j+1} - t_j) \leq h$, the unique solution of the problem (2.4), (2.5), (2.10), (2.20), (2.16) with initial condition given by (2.9) and $w[0] = w_0$ satisfies the following estimate for $t \geq 0$ for the observer error $e[t] = w[t] - u[t]$:*

$$\|e[t]\| \leq (1-\Omega)^{-1}\sqrt{\max\left(|P|, \dfrac{Q}{2}\right)}\exp(-\kappa t)\|e[0]\|$$

$$+ \exp(\kappa h)(1-\Omega)^{-1}\gamma \sum_{i=1}^{m}\left(\|l_i\| + h\sum_{r=1}^{m}\|l_r\|\left|\int_0^1 c_r(x)l_i(x)dx\right|\right)\sup_{0 \leq s \leq t}\left(|\xi_i(s)|\exp(-\kappa(t-s))\right) \quad (2.25)$$

$$+ (1-\Omega)^{-1}\gamma\left(1 + h\exp(\kappa h)\sum_{i=1}^{m}\|l_i\|\|c_i\|\right)\sup_{0 \leq s \leq t}\left(\|v[s] - \tilde{v}[s]\|\exp(-\kappa(t-s))\right)$$

**Remark 2.4:** **(a)** The proofs of Theorem 2.2 and Theorem 2.3 are provided in Section 4; they are based on a small-gain argument as well as Lyapunov analysis. More specifically, by using an appropriate Lyapunov functional for an auxiliary problem, we are in a position to utilize the Input-to-Output Stability (IOS) property and prove the desired estimates (2.23) and (2.25).
**(b)** It should be noticed that estimates (2.23), (2.25) are IOS estimate for the observer error with respect to measurement noise and modeling errors. When measurement noise and modeling errors are absent, estimates (2.23), (2.25) imply global exponential convergence of the observer error.
**(c)** The constant $h > 0$ for which the inequality $\sup_{j \geq 0}(t_{j+1} - t_j) \leq h$ holds as well as (2.22) or (2.24) is the diameter of the sampling schedule. In general, inequalities (2.22), (2.24) impose bounds on the diameter of the sampling schedule. Moreover, it should be noticed that inequality (2.24) is more demanding than inequality (2.22) in the sense that all $h > 0$ and $\kappa \in [0, \mu)$ that satisfy (2.24) necessarily also satisfy automatically (2.22). That is why the sampled-data observer with the inter-sample predictor allows less frequent sampling than the sampled-data observer without the inter-sample predictor.
**(d)** Inequalities (2.23), (2.25) show that the sampled-data observer without the inter-sample predictor is more sensitive to measurement noise than the sampled-data observer with the inter-sample predictor.
**(e)** For both Theorem 2.2 and Theorem 2.3, the existence of a symmetric, positive definite matrix $P \in \Re^{N \times N}$ with $P \geq I$ for which there exists a constant $\sigma > 0$ such that $PA + A'P \leq -2\sigma P$ is not an issue, since it is assumed that the real matrix $A \in \Re^{N \times N}$ defined by (2.12) is a Hurwitz matrix.

## 3. Illustrative Examples

The examples presented in this section have multiple purposes. First of all, the examples aim to show how easily we can apply Theorem 2.2 in order to design sampled-data observers for 1-D parabolic systems. Furthermore, the examples also illustrate the following additional facts:
- the first example shows that there exist parabolic PDEs which allow the diameter of the sampling schedule to be arbitrarily large for the observer with the inter-sample predictor, while the observer without the inter-sample predictor requires a sufficiently frequent sampling schedule, and
- the second example shows that Theorem 2.2 can be used also when the outputs are not non-local outputs of the form (2.10) and can guarantee stronger estimates of the observer error (e.g., estimates in the spatial sup norm).

**Example 3.1:** Consider the PDE problem

$$\frac{\partial u}{\partial t}(t,x) = p\frac{\partial^2 u}{\partial x^2}(t,x) + v(t,x), \text{ for } (t,x) \in (0,+\infty) \times (0,1) \quad (3.1)$$

where $p > 0$ is a constant and $v: \Re_+ \times [0,1] \to \Re$ is an input, under Neumann boundary conditions

$$\frac{\partial u}{\partial x}(t,0) = \frac{\partial u}{\partial x}(t,1) = 0, \text{ for } t \geq 0 \quad (3.2)$$

and a sampled, scalar, non-local output

$$y(t) = \xi(t_j) + \int_0^1 xu(t_j,x)dx, \text{ for } t \in [t_j, t_{j+1}) \text{ and } j = 0,1,2,... \quad (3.3)$$

where $\xi: \Re_+ \to \Re$ is the measurement error and $\{t_j \geq 0 : j = 0,1,2,...\}$ (the sampling schedule) is an increasing sequence of times (the sampling times) with $t_0 = 0$ and $\lim_{j \to +\infty}(t_j) = +\infty$. It is clear that system (3.1), (3.2), (3.3) is a system of the form (2.4), (2.5), (2.10) with $q(x) \equiv 0$, $K(u) \equiv 0$, $g(x, \bar{P}u) \equiv 0$, $a_0 = a_1 = 0$, $b_0 = b_1 = 1$, $m = 1$ and $k_1(x) = x$.

The eigenvalues and eigenfunctions of the SL operator $B: D \to L^2(0,1)$ defined by (2.1), (2.2), are

$$\lambda_1 = 0 \text{ and } \lambda_n = p(n-1)^2\pi^2 \text{ for } n \geq 2 \quad (3.4)$$

$$\phi_1(x) \equiv 1 \text{ and } \phi_n(x) = \sqrt{2}\cos(n\pi x) \text{ for } n \geq 2 \text{ and } x \in [0,1] \quad (3.5)$$

It follows that Assumptions (H1), (H2) hold. More specifically, inequality (2.7) holds with $R = 0$. We next apply Theorem 2.2 with $N = 1$, $c_1(x) = 1/2$, $P = [1]$ and $L_{1,1} = -p\pi^2$. Definition (2.13) in conjunction with (3.5) implies that $c_{1,j} = 0$, for $j \geq 2$ and $c_{1,1} = 1/2$. Definitions (2.12), (2.14) in conjunction with (3.4), (3.5), give $A_{1,1} = -p\pi^2/2$, $l_i(x) \equiv -p\pi^2$. The proposed sampled-data observer with the inter-sample predictor takes the form

$$\frac{\partial w}{\partial t}(t,x) = p\frac{\partial^2 w}{\partial x^2}(t,x) + \tilde{v}(t,x) - p\pi^2\left(\frac{1}{2}\int_0^1 w(t,s)ds - \zeta(t)\right), \text{ for } (t,x) \in (0,+\infty) \times (0,1) \quad (3.6)$$

$$\frac{\partial w}{\partial x}(t,0) = \frac{\partial w}{\partial x}(t,1) = 0, \text{ for } t \geq 0 \quad (3.7)$$

$$\zeta(t_j) = y(t_j) - \int_0^1\left(x - \frac{1}{2}\right)w(t_j,x)dx, \text{ for } j = 0,1,2,... \quad (3.8)$$

$$\dot{\zeta}(t) = \frac{1}{2}\int_0^1 \tilde{v}(t,x)dx, \text{ for } t \in [t_j, t_{j+1}) \text{ and } j = 0,1,2,... \quad (3.9)$$

Using all the above, we conclude that all assumptions of Theorem 2.2 hold with $R = 0$, $K = 0$, $\sigma = \mu = p\pi^2/2$, $Q = 2$, $H(Q) \equiv 0$, $\|k_1 - c_1\| = \frac{1}{2\sqrt{3}}$, $\kappa = \mu\omega$, $\tilde{g} := \frac{2}{p\pi^2}$, $\gamma := \frac{1}{p\pi^2}\sqrt{\frac{2}{1-\omega}}$, for all $h > 0$ and $\omega \in [0,1)$ for which the small-gain condition

$$\Omega := \frac{\exp\left(\omega\frac{p}{2}\pi^2 h\right)}{\sqrt{6(1-\omega)}} < 1 \quad (3.10)$$

holds. Therefore, for all $h > 0$ and $\omega \in [0,1)$ for which (3.10) holds, the following property also holds: for every $u_0, w_0 \in \left\{\theta \in H^2(0,1) : \frac{d\theta}{dx}(0) = \frac{d\theta}{dx}(1) = 0\right\}$, for every bounded mapping $\xi: \Re_+ \to \Re$, for every pair of inputs $v, \tilde{v} \in C^0(\Re_+ \times [0,1])$ that satisfy Assumption (H3) and for every increasing sequence of times $\{t_j \geq 0: j=0,1,2,...\}$ with $t_0 = 0$, $\lim_{j \to +\infty}(t_j) = +\infty$ and $\sup_{j \geq 0}(t_{j+1} - t_j) \leq h$, the unique solution of the initial-boundary value problem (3.1), (3.2), (3.3), (3.6), (3.7), (3.8), (3.9) with initial condition given by (2.9) and $w[0] = w_0$ satisfies the estimate for all $t \geq 0$

$$\left(\sqrt{6(1-\omega)} - \exp\left(\frac{p\omega}{2}\pi^2 h\right)\right)\|e[t]\| \leq \sqrt{6(1-\omega)} \exp\left(-\frac{p\omega}{2}\pi^2 t\right)\|e[0]\|$$

$$+ \sqrt{12} \exp\left(\frac{p\omega}{2}\pi^2 h\right) \sup_{0 \leq s \leq t}\left(|\xi(s)|\exp\left(-\frac{p\omega}{2}\pi^2(t-s)\right)\right) \quad (3.11)$$

$$+ \sqrt{12}\left(\frac{1}{p\pi^2} + \frac{h}{2}\exp\left(\frac{p\omega}{2}\pi^2 h\right)\right) \sup_{0 \leq s \leq t}\left(\|v[s] - \tilde{v}[s]\|\exp\left(-\frac{p\omega}{2}\pi^2(t-s)\right)\right)$$

where $e[t] = w[t] - u[t]$. It should be noticed that for every given $h > 0$ there exists $\omega^* \in (0,1)$ such that (3.10) holds for all $\omega \in [0, \omega^*)$. Therefore, for every sampling schedule the observer error will converge to zero in absence of noise and modeling errors. However, notice that a large value for $h > 0$ (i.e., when measurements are sparse) will give a small value for $\omega^* \in (0,1)$, i.e., a slow convergence of the observer error. Moreover, (3.11) shows that a large value for $h > 0$ implies sensitivity with respect to modeling errors since the gain coefficient of $v - \tilde{v}$ increases with $h > 0$.

On the other hand, the observer without the inter-sample predictor is given by (3.7) with

$$\frac{\partial w}{\partial t}(t,x) = p\frac{\partial^2 w}{\partial x^2}(t,x) + \tilde{v}(t,x) - p\pi^2\left(\int_0^1 sw(t_j,s)ds - y(t_j)\right), \text{ for } (t,x) \in (t_j, t_{j+1}) \times (0,1) \text{ and } j = 0,1,2,... \quad (3.12)$$

Using all the above, we conclude that all assumptions of Theorem 2.3 hold with $N=1$, $c_1(x) = 1/2$, $P = [1]$, $L_{1,1} = -p\pi^2$, $A_{1,1} = -p\pi^2/2$, $l_i(x) \equiv -p\pi^2$, $R = 0$, $K = 0$, $\sigma = \mu = p\pi^2/2$, $Q = 2$, $H(Q) \equiv 0$, $\|k_1 - c_1\| = \frac{1}{2\sqrt{3}}$, $\kappa = \mu\omega$, $\tilde{g} := \frac{2}{p\pi^2}$, $\gamma := \frac{1}{p\pi^2}\sqrt{\frac{2}{1-\omega}}$, for all $h > 0$ and $\omega \in [0,1)$ for which the condition

$$\Omega := \exp\left(\omega\frac{p\pi^2}{2}\right)\frac{hp\pi^2 + 1}{\sqrt{6(1-\omega)}} < 1 \quad (3.13)$$

holds. It is clear that in this case the upper diameter of the sampling schedule is not allowed to be greater or equal to $(\sqrt{6} - 1)/(p\pi^2)$. For uniform sampling schedules ($t_j = jh$ for $j = 0,1,2,...$) we are in a position to give necessary and sufficient conditions for the successful operation of the sampled-data observer (3.7), (3.12): the sampling period $h > 0$ has to be strictly less than $4/(p\pi^2)$. Therefore, the fact that the sampled-data observer (3.7), (3.12) requires a sufficiently small upper diameter of the sampling schedule is not an artifact of the analysis but is a fundamental limitation of the observer (3.7), (3.12). Thus, as stated in Remark 2.4(c) the observer without the inter-sample predictor requires a sufficiently frequent sampling schedule (here in sharp contrast with the sampled-data observer with the inter-sample predictor). ◁

**Example 3.2:** Consider the PDE problem

$$\frac{\partial u}{\partial t}(t,x) = p\frac{\partial^2 u}{\partial x^2}(t,x) - qu(t,x) + v(t,x), \text{ for } (t,x) \in (0,+\infty) \times (0,1) \quad (3.14)$$

where $p > 0$, $q \in \Re$ are constants, $v: \Re_+ \times [0,1] \to \Re$ is an input, under the following boundary conditions

$$u(t,0) = \frac{\partial u}{\partial x}(t,1) = 0, \text{ for } t \geq 0 \quad (3.15)$$

and a sampled, scalar, local output

$$y(t) = \xi(t_j) + u(t_j, 1), \text{ for } t \in [t_j, t_{j+1}) \text{ and } j = 0,1,2,... \quad (3.16)$$

where $\xi: \Re_+ \to \Re$ is the measurement error and $\{t_j \geq 0: j = 0,1,2,...\}$ (the sampling schedule) is an increasing sequence of times (the sampling times) with $t_0 = 0$ and $\lim_{j \to +\infty}(t_j) = +\infty$. The output gives the boundary point value of the state and it is not a non-local output of the form (2.10). System (3.14), (3.15), (3.16) is also studied in [30] (with $x$ replaced by $1-x$).

Despite the fact that the output is not a non-local output of the form (2.10), we show next that Theorem 2.2 can be used for the design of a sampled-data observer. To make things simpler, we assume that the reaction coefficient $q \in \Re$ satisfies the inequality

$$-9p\pi^2 < 4q < 7p\pi^2 \quad (3.17)$$

although the observer can be designed even when the reaction coefficient does not satisfy (3.17).

In order to be able to apply Theorem 2.2 with a local measurement, we need to look at the variable

$$\tilde{u}(t,x) = \frac{\partial u}{\partial x}(t,x) + p(x-1)v(t,0), \text{ for } t \geq 0, \ x \in [0,1] \quad (3.18)$$

which contains information for the spatial derivative of the state $u$ and not the state itself. Indeed, Proposition 5.11 on page 113 in [19] guarantees that (under sufficient regularity for the initial condition of the state and the input) $\tilde{u}$ defined by (3.18), satisfies the PDE

$$\frac{\partial \tilde{u}}{\partial t}(t,x) = p\frac{\partial^2 \tilde{u}}{\partial x^2}(t,x) - q\tilde{u}(t,x) + \tilde{v}(t,x), \text{ for } (t,x) \in (0,+\infty) \times (0,1) \quad (3.19)$$

and the boundary conditions

$$\frac{\partial \tilde{u}}{\partial x}(t,0) = \tilde{u}(t,1) = 0, \text{ for } t \geq 0 \quad (3.20)$$

where

$$\tilde{v}(t,x) := pq(x-1)v(t,0) + \frac{\partial v}{\partial x}(t,x) + p(x-1)\frac{\partial v}{\partial t}(t,0), \text{ for } t \geq 0, \ x \in [0,1] \quad (3.21)$$

Moreover, the output map (3.16) can be expressed (using (3.18)) in the following way:

$$y(t) = \xi(t_j) + \int_0^1 \tilde{u}(t_j,x)dx + \frac{p}{2}v(t_j,0), \text{ for } t \in [t_j, t_{j+1}) \text{ and } j = 0,1,2,... \quad (3.22)$$

It is exactly this formulation to which Theorem 2.2 can be applied: system (3.19), (3.20), (3.22) is a system of the form (2.4), (2.5), (2.10) with $q(x) \equiv q$, $K(u) \equiv 0$, $g(x, \bar{P}u) \equiv 0$, $a_0 = b_1 = 0$, $b_0 = a_1 = 1$, $m = 1$ and $k_1(x) \equiv 1$. The eigenvalues and eigenfunctions of operator $B$ defined by (2.1), (2.2), are

$$\lambda_n = p\frac{(2n-1)^2}{4}\pi^2 + q \text{ for } n \geq 1 \quad (3.23)$$

$$\phi_n(x) = \sqrt{2}\cos\left(\frac{2n-1}{2}\pi x\right) \text{ for } n \geq 1 \text{ and } x \in [0,1] \quad (3.24)$$

It follows that Assumptions (H1), (H2) hold. More specifically, inequality (2.7) holds with $R = 0$. We next apply Theorem 2.2 with $N = 1$, $c_1(x) = \frac{4}{\pi}\cos\left(\frac{\pi x}{2}\right)$, $P = [1]$ and $L_{1,1} = \pi\frac{4q - 7p\pi^2}{16\sqrt{2}}$. Definition (2.13) in conjunction with (3.24) implies that $c_{1,j} = 0$, for $j \geq 2$ and $c_{1,1} = 2\sqrt{2}/\pi$. Definitions (2.12), (2.14) in conjunction with (3.23), (3.24), give $A_{1,1} = -p\frac{9\pi^2}{8} - \frac{q}{2}$, $l_1(x) = \pi\frac{4q - 7p\pi^2}{16}\cos\left(\frac{\pi x}{2}\right)$. Notice that by virtue of (3.17), it holds that $A_{1,1} < 0$. The proposed sampled-data observer takes the form

$$\frac{\partial w}{\partial t}(t,x) = p\frac{\partial^2 w}{\partial x^2}(t,x) - qw(t,x) + \tilde{v}(t,x) + \pi\frac{4q - 7p\pi^2}{16}\cos\left(\frac{\pi x}{2}\right)\left(\frac{4}{\pi}\int_0^1 \cos\left(\frac{\pi s}{2}\right)w(t,s)ds - \zeta_1(t)\right)$$

$$\text{for } (t,x) \in (0,+\infty) \times (0,1) \quad (3.25)$$

$$\frac{\partial w}{\partial x}(t,0) = w(t,1) = 0, \text{ for } t \geq 0 \tag{3.26}$$

$$\zeta_1(t_j) = y(t_j) - \frac{p}{2}v(t_j,0) - \int_0^1 \left(1 - \frac{4}{\pi}\cos\left(\frac{\pi x}{2}\right)\right)w(t_j,x)dx, \text{ for } j = 0,1,2,... \tag{3.27}$$

$$\dot{\zeta}_1(t) = -\frac{4}{\pi}\int_0^1 \left(\left(p\frac{\pi^2}{4} + q\right)w(t,x) - \tilde{v}(t,x)\right)\cos\left(\frac{\pi x}{2}\right)dx, \text{ for } t \in [t_j, t_{j+1}) \text{ and } j = 0,1,2,... \tag{3.28}$$

However, since the observer state $w$ will not estimate the state $u$ but the variable $\tilde{u}$ (defined by (3.18)), we also need to get an estimation of the state $u$. This will be achieved by the state estimator

$$\hat{u}(t,x) = \int_0^x w(t,s)ds - pv(t,0)\left(\frac{x^2}{2} - x\right), \text{ for } t \geq 0, \ x \in [0,1] \tag{3.29}$$

Using all the above, we conclude that all assumptions of Theorem 2.2 hold with $K = 0$, $\sigma = \mu = p\frac{9\pi^2}{8} + \frac{q}{2}$, $Q = 2$, $H(Q) \equiv 0$, $\|k_1 - c_1\| = \frac{\sqrt{2}}{\pi}\sqrt{\pi^2 - 8}$, $\kappa = \mu\omega$, $\tilde{g} := \frac{8}{9p\pi^2 + 4q}$, $\gamma := \frac{4\sqrt{2}}{(9p\pi^2 + 4q)\sqrt{1-\omega}}$, for all $h > 0$ and $\omega \in [0,1)$ for which the small-gain condition

$$\Omega(\omega,h) = \exp\left(\omega h \frac{9p\pi^2 + 4q}{8}\right)\frac{(7p\pi^2 - 4q)}{2\sqrt{2}(9p\pi^2 + 4q)\sqrt{1-\omega}}\left(\frac{|p\pi^2 + 4q|}{4\sqrt{2}}\pi h + \sqrt{\pi^2 - 8}\right) < 1 \tag{3.30}$$

holds. By virtue of continuity of $\Omega(\omega,h)$ at $(\omega,h) = (0,0)$ and the fact that $\Omega(0,0) < 1$, it follows that there exist $h > 0$ and $\omega \in [0,1)$ for which (3.30) holds. Therefore, for all $h > 0$ and $\omega \in [0,1)$ for which (3.30) holds, the following property also holds: there exists a constant $\Theta > 0$ such that for every $\tilde{u}_0, w_0 \in \left\{\theta \in H^2(0,1) : \frac{d\theta}{dx}(0) = \theta(1) = 0\right\}$, for every bounded mapping $\xi : \Re_+ \to \Re$, for every input $v \in C^1(\Re_+ \times [0,1])$ for which $\tilde{v}$ being defined by (3.21) satisfies Assumption (H3) and for every increasing sequence of times $\{t_j \geq 0 : j = 0,1,2,...\}$ with $t_0 = 0$, $\lim_{j \to +\infty}(t_j) = +\infty$ and $\sup_{j \geq 0}(t_{j+1} - t_j) \leq h$, the unique solution of the initial-boundary value problem (3.19), (3.20), (3.22), (3.25), (3.26), (3.27), (3.28) with initial condition $\tilde{u}[0] = \tilde{u}_0$, $w[0] = w_0$ satisfies the estimate

$$\|e[t]\| \leq \Theta \exp(-\kappa t)\|e[0]\| + \Theta \sup_{0 \leq s \leq t}(|\xi(s)|), \text{ for all } t \geq 0 \tag{3.31}$$

where $e[t] = w[t] - \tilde{u}[t]$. Estimate (3.31) combined with (3.15), (3.18) and (3.29) implies the estimate

$$\max_{0 \leq x \leq 1}(|\hat{u}(t,x) - u(t,x)|) \leq \Theta \exp(-\kappa t)\|e[0]\| + \Theta \sup_{0 \leq s \leq t}(|\xi(s)|), \text{ for all } t \geq 0 \tag{3.32}$$

which is an estimate of the observer error in the spatial sup norm rather than the $L^2$ norm. Estimates (3.31), (3.32) should be compared with the estimate $L^2$ norm observer error estimate that is provided in [30] (depending also on the $H^1$ norm of the initial error). ◁

## 4. Proofs of Main Results

For the proofs of the main results we need the following auxiliary result.

**Proposition 4.1:** *Consider the SL operator $B$ under Assumption (H1). Let $N, m \geq 1$ be integers with $\lambda_{N+1} > 0$, $c_i \in D$ ($i = 1,...,m$) be given functions and let $L = \{L_{i,j} : i = 1,...,N, j = 1,...,m\} \in \Re^{N \times m}$ be a matrix so that the matrix $A \in \Re^{N \times N}$ defined by (2.12) is a Hurwitz matrix. Define $l_i \in D$ for $i = 1,...,m$ by means of (2.14) and the constant $K$ by means of (2.21). Let $P \in \Re^{N \times N}$ be a symmetric, positive* definite matrix with $P \geq I$ for which there exists a constant $\sigma > 0$ such that $PA + A'P \leq -2\sigma P$. Let $e_0 \in D$, $T > 0$ and let $\bar{v} \in C^0([0,T] \times [0,1])$ be an input. Then every solution $e \in C^0([0,T] \times [0,1]) \cap C^1((0,T] \times [0,1])$ with $e[t] \in C^2([0,1]) \cap D$ for $t \in (0,T]$ of the problem

$$\frac{\partial e}{\partial t}(t,x) = p\frac{\partial^2 e}{\partial x^2}(t,x) - q(x)e(t,x) + \sum_{i=1}^{m} l_i(x)\int_0^1 c_i(s)e(t,s)ds + \bar{v}(t,x), \text{ for } t \in (0,T], \ x \in (0,1) \quad (4.1)$$

$$a_0 e(t,0) + b_0 \frac{\partial e}{\partial x}(t,0) = a_1 e(t,1) + b_1 \frac{\partial e}{\partial x}(t,1) = 0, \text{ for } t \in [0,T] \quad (4.2)$$

$$e[0] = e_0 \quad (4.3)$$

satisfies the following estimate for every $Q \geq 2$ with $Q > 2|L'PL|K^2/(\sigma \lambda_{N+1})$:

$$\|e[t]\|\exp(\kappa t) \leq \sqrt{\max\left(|P|, \frac{Q}{2}\right)}\|e[0]\| + \sqrt{\frac{\tilde{g}}{2(\mu - \kappa)}}\sup_{0 \leq s \leq t}\left(\|\bar{v}[s]\|\exp(\kappa s)\right), \text{ for all } t \in [0,T], \ \kappa \in [0,\mu) \quad (4.4)$$

where $\tilde{g} := \max\left(\frac{4|P|}{4\sigma + H(Q)}, \frac{Q}{2\lambda_{N+1}}\right)$, $\mu := \frac{H(Q) + 2\lambda_{N+1}}{4}$ and $H(Q) := 2\sigma - \lambda_{N+1} - \sqrt{(2\sigma - \lambda_{N+1})^2 + 16Q^{-1}|L'PL|K^2}$.

**Proof:** Let $Q \geq 2$ with $Q > 2|L'PL|K^2/(\sigma \lambda_{N+1})$ and let $e \in C^0([0,T] \times [0,1]) \cap C^1((0,T] \times [0,1])$ with $e[t] \in C^2([0,1]) \cap D$ for $t \in (0,T]$ be a solution of the initial-boundary value problem (4.1), (4.2), (4.3). Define the Lyapunov functional $V : [0,T] \to \Re_+$ by the formula:

$$V(t) := \xi'(t)P\xi(t) + \frac{Q}{2}\sum_{n=N+1}^{\infty} r_n^2(t) \quad (4.5)$$

where

$$r_n(t) = \int_0^1 e(t,x)\phi_n(x)dx, \text{ for } n = 1,2,... \text{ and } \xi(t) = (r_1(t),...,r_N(t))' \quad (4.6)$$

It should be noticed that by virtue of (4.6) and (4.1), (4.2) the following equations hold for $t \in (0,T]$:

$$\dot{r}_n(t) = -\lambda_n r_n(t) + L_{n,1}\int_0^1 c_1(x)e(t,x)dx + ... + L_{n,m}\int_0^1 c_m(x)e(t,x)dx + v_n(t), \text{ for } n = 1,...,N \quad (4.7)$$

$$\dot{r}_n(t) = -\lambda_n r_n(t) + v_n(t), \text{ for } n = N+1,... \quad (4.8)$$

where

$$v_n(t) = \int_0^1 \bar{v}(t,x)\phi_n(x)dx, \text{ for } n = 1,2,... \quad (4.9)$$

Using (4.6), (2.12) and the fact that $\int_0^1 c_i(x)e(t,x)dx = \sum_{j=1}^{\infty} c_{i,j}r_j(t)$ for $i = 1,...,m$ (a direct consequence of (2.13) and (4.6)), we get that the following differential equations hold for $t \in (0,T]$:

$$\dot{\xi}(t) = A\xi(t) + F(t) \quad (4.10)$$

where

$$F(t) = \begin{bmatrix} v_1(t) \\ \vdots \\ v_N(t) \end{bmatrix} + L\begin{bmatrix} \sum_{j=N+1}^{\infty} c_{1,j}r_j(t) \\ \vdots \\ \sum_{j=N+1}^{\infty} c_{m,j}r_j(t) \end{bmatrix} \quad (4.11)$$

Using (4.8), (4.10), definition (4.5) and the fact that $PA + A'P \leq -2\sigma P$, we obtain for $t \in (0,T]$:

$$\dot{V}(t) \leq -2\sigma\xi'(t)P\xi(t) + 2\xi'(t)PF(t) - Q\sum_{n=N+1}^{\infty} \lambda_n r_n^2(t) + Q\sum_{n=N+1}^{\infty} r_n(t)v_n(t) \quad (4.12)$$

Using the fact that $\lambda_n > 0$ for $n \geq N+1$ (since $\lambda_{N+1} > 0$) and the fact that $2r_n(t)v_n(t) \leq \lambda_n r_n^2(t) + \lambda_n^{-1} r_n^2(t)$ for $n \geq N+1$, we obtain from (4.12) for $t \in (0,T]$:

$$\dot{V}(t) \leq -2\sigma\xi'(t)P\xi(t) + 2\xi'(t)PF(t) - \frac{Q}{2}\sum_{n=N+1}^{\infty}\lambda_n r_n^2(t) + Q\sum_{n=N+1}^{\infty}\frac{1}{2\lambda_n}v_n^2(t) \qquad (4.13)$$

The fact that $\lambda_n \geq \lambda_{N+1}$ for $n \geq N+1$ in conjunction with (4.13) and (4.11) gives for $t \in (0,T]$:

$$\dot{V}(t) \leq -2\sigma\xi'(t)P\xi(t) + 2\xi'(t)PF(t) - \frac{Q}{2}\lambda_{N+1}\sum_{n=N+1}^{\infty} r_n^2(t) + \frac{Q}{2\lambda_{N+1}}\sum_{n=N+1}^{\infty} v_n^2(t) \qquad (4.14)$$

For any two vectors $x, y \in \Re^N$ and any $\varepsilon > 0$, it holds that $2x'y \leq \varepsilon x'P^{-1}x + \varepsilon^{-1}y'Py$. Applying this fact we get from (4.14) for every $t \in (0,T]$ and $\varepsilon > 0$:

$$\dot{V}(t) \leq -2(\sigma - \varepsilon)\xi'(t)P\xi(t) + \varepsilon^{-1}\begin{bmatrix} v_1(t) \\ \vdots \\ v_N(t) \end{bmatrix}' P \begin{bmatrix} v_1(t) \\ \vdots \\ v_N(t) \end{bmatrix} + \varepsilon^{-1}\begin{bmatrix} \sum_{j=N+1}^{\infty} c_{1,j}r_j(t) \\ \vdots \\ \sum_{j=N+1}^{\infty} c_{m,j}r_j(t) \end{bmatrix}' L'PL \begin{bmatrix} \sum_{j=N+1}^{\infty} c_{1,j}r_j(t) \\ \vdots \\ \sum_{j=N+1}^{\infty} c_{m,j}r_j(t) \end{bmatrix} \qquad (4.15)$$

$$-\frac{Q}{2}\lambda_{N+1}\sum_{n=N+1}^{\infty} r_n^2(t) + \frac{Q}{2\lambda_{N+1}}\sum_{n=N+1}^{\infty} v_n^2(t)$$

Using the fact that $\left|\sum_{j=N+1}^{\infty} c_{i,j}r_j(t)\right| \leq \left(\sum_{j=N+1}^{\infty} c_{i,j}^2\right)^{1/2}\left(\sum_{j=N+1}^{\infty} r_j^2(t)\right)^{1/2}$ for $i=1,...,m$ in conjunction with definition (2.21), we get from (4.15) for every $t \in (0,T]$ and $\varepsilon > 0$:

$$\dot{V}(t) \leq -2(\sigma - \varepsilon)\xi'(t)P\xi(t) - 2\left(\frac{\lambda_{N+1}}{2} - \varepsilon^{-1}Q^{-1}|L'PL|K^2\right)\frac{Q}{2}\sum_{n=N+1}^{\infty} r_n^2(t) + \varepsilon^{-1}|P|\sum_{n=1}^{N} v_n^2(t) + \frac{Q}{2\lambda_{N+1}}\sum_{n=N+1}^{\infty} v_n^2(t) \qquad (4.16)$$

Setting $\varepsilon = \left(2\sigma - \lambda_{N+1} + \sqrt{(2\sigma - \lambda_{N+1})^2 + 16Q^{-1}|L'PL|K^2}\right)/4$, we get from (4.16) and (4.5) for $t \in (0,T]$:

$$\dot{V}(t) \leq -2\mu V(t) + \tilde{g}\sum_{n=1}^{\infty} v_n^2(t) \qquad (4.17)$$

where $\tilde{g} := \max\left(\frac{4|P|}{4\sigma + H(Q)}, \frac{Q}{2\lambda_{N+1}}\right)$, $\mu := \frac{H(Q) + 2\lambda_{N+1}}{4}$ and $H(Q) := 2\sigma - \lambda_{N+1} - \sqrt{(2\sigma - \lambda_{N+1})^2 + 16Q^{-1}|L'PL|K^2}$.

Notice that Parseval's identity and (4.9) give $\|\bar{v}[t]\|^2 = \sum_{n=1}^{\infty} v_n^2(t)$ for $t \geq 0$. Therefore (4.17) in conjunction with continuity of the mapping $t \to V(t)$ for $t \geq 0$ implies the following inequality for $t \in [0,T]$:

$$V(t) \leq \exp(-2\mu t)V(0) + \tilde{g}\int_0^t \exp(-2\mu(t-s))\|\bar{v}[s]\|^2 ds \qquad (4.18)$$

Notice that Parseval's identity in conjunction with (4.6) implies that $\|e[t]\|^2 = \sum_{n=1}^{\infty} r_n^2(t)$. Consequently, the facts that $P \geq I$, $Q \geq 2$ in conjunction with (4.5) imply the following inequalities for $t \in [0,T]$:

$$\|e[t]\|^2 \leq V(t) \text{ and } V(0) \leq \max\left(|P|, \frac{Q}{2}\right)\|e[0]\|^2 \qquad (4.19)$$

Combining (4.18), (4.19), we get for all $t \in [0,T]$, $\kappa \in [0,\mu)$:

$$2\|e[t]\|^2 \exp(2\kappa t) \leq \exp(-2(\mu - \kappa)t)\max(2|P|, Q)\|e[0]\|^2 + \tilde{g}\sup_{0 \leq s \leq t}\left(\|\bar{v}[s]\|^2 \exp(2\kappa s)\right)/(\mu - \kappa) \qquad (4.20)$$

Estimate (4.4) is a direct consequence of estimate (4.20). The proof is complete. ◁

We are now ready to give the proofs of the main results of the present work.

**Proof of Theorem 2.2:** Let $u_0, w_0 \in D$, a bounded mapping $\xi : \Re_+ \to \Re^m$, inputs $v, \tilde{v} \in C^0(\Re_+ \times [0,1])$ satisfying (H3) and an increasing sequence $\{t_j \geq 0 : j = 0,1,2,...\}$ with $t_0 = 0$, $\lim_{j \to +\infty}(t_j) = +\infty$ and $\sup_{j \geq 0}(t_{j+1} - t_j) \leq h$, where $h > 0$ satisfies (2.22). Consider the solution $u[t], w[t]$, $\zeta_i(t)$ ($i = 1,...,m$) of (2.4), (2.5), (2.10), (2.15), (2.16), (2.18), (2.19) with initial condition (2.9) and $w[0] = w_0$. Define for $t \geq 0$:

$$\varepsilon_i(t) := \zeta_i(t) - \int_0^1 c_i(x) u(t,x) dx \quad (i=1,...,m). \tag{4.21}$$

$$e[t] = w[t] - u[t] \tag{4.22}$$

Moreover, define for all $0 \leq a \leq b$:

$$\|e\|_{[a,b]} := \sup_{a \leq s \leq b} (\|e[s]\| \exp(\kappa s)), \quad \|v - \tilde{v}\|_{[a,b]} := \sup_{a \leq s \leq b} (\|v[s] - \tilde{v}[s]\| \exp(\kappa s))$$
$$\|\varepsilon_i\|_{[a,b]} := \sup_{a \leq s \leq b} (|\varepsilon_i(s)| \exp(\kappa s)), \quad \|\xi_i\|_{[a,b]} := \sup_{a \leq s \leq b} (|\xi_i(s)| \exp(\kappa s)), \quad i = 1,...,m \tag{4.23}$$

Using (2.17), (2.19) and (4.21), (4.22) we get for all $t \in \overline{I} := \Re_+ \setminus \{t_j : j = 0,1,2,...\}$, $i = 1,...,m$:

$$\dot{\varepsilon}_i(t) = \int_0^1 \left( p \frac{d^2 c_i}{dx^2}(x) - q(x) c_i(x) \right) e(t,x) dx + \int_0^1 c_i(x)(f(w[t]) - f(u[t]))(x) dx + \int_0^1 c_i(x)(\tilde{v}(t,x) - v(t,x)) dx \tag{4.24}$$

Let an arbitrary integer $j \geq 0$ and $t \in (t_j, t_{j+1})$. It follows from (4.24), (2.7) and the Cauchy-Schwarz inequality that the following estimate holds for all $\tau \in (t_j, t]$:

$$|\dot{\varepsilon}_i(\tau)| \leq \left( \left\| p \frac{d^2 c_i}{dx^2} - q c_i \right\| + R \|c_i\| \right) \max_{t_j \leq s \leq t}(\|e[s]\|) + \|c_i\| \max_{t_j \leq s \leq t}(\|v[s] - \tilde{v}[s]\|) \tag{4.25}$$

Therefore, we obtain from (4.25) the estimate:

$$|\varepsilon_i(t)| \leq |\varepsilon_i(t_j)| + \left( \left\| p \frac{d^2 c_i}{dx^2} - q c_i \right\| + R \|c_i\| \right)(t - t_j) \max_{t_j \leq s \leq t}(\|e[s]\|) + \|c_i\|(t - t_j) \max_{t_j \leq s \leq t}(\|v[s] - \tilde{v}[s]\|) \tag{4.26}$$

Using (2.18), (2.10), (4.21) and (4.22) we get $\varepsilon_i(t_j) = \xi_i(t_j) - \int_0^1 (k_i(x) - c_i(x)) e(t_j, x) dx$, for $i = 1,...,m$, which (by means of the Cauchy-Schwarz inequality) gives

$$|\varepsilon_i(t_j)| \leq |\xi_i(t_j)| + \|k_i - c_i\| \|e[t_j]\|, \text{ for } i = 1,...,m \tag{4.27}$$

Combining (4.26), (4.27) we obtain the following estimate:

$$|\varepsilon_i(t)| \exp(\kappa t) \leq |\xi_i(t_j)| \exp(\kappa t) + \|c_i\|(t - t_j) \exp(\kappa(t - t_j)) \sup_{t_j \leq s \leq t}(\|v[s] - \tilde{v}[s]\| \exp(\kappa s))$$
$$+ \|k_i - c_i\| \|e[t_j]\| \exp(\kappa t) + \left( \left\| p \frac{d^2 c_i}{dx^2} - q c_i \right\| + R \|c_i\| \right)(t - t_j) \exp(\kappa(t - t_j)) \sup_{t_j \leq s \leq t}(\|e[s]\| \exp(\kappa s)) \tag{4.28}$$

Using the facts that $t \in [t_j, t_{j+1})$, $\sup_{j \geq 0}(t_{j+1} - t_j) \leq h$ and (4.28), (4.23), we get for $t \geq 0$:

$$\|\varepsilon_i\|_{[0,t]} \leq \exp(\kappa h) \|\xi_i\|_{[0,t]} + \exp(\kappa h) \left( \left( \left\| p \frac{d^2 c_i}{dx^2} - q c_i \right\| + R \|c_i\| \right) h + \|k_i - c_i\| \right) \|e\|_{[0,t]} + h \exp(\kappa h) \|c_i\| \|v - \tilde{v}\|_{[0,t]} \tag{4.29}$$

It follows from (2.4), (2.5), (2.6), (2.20), (2.16), (4.21), (4.22) that $e[t]$ is a solution of (4.1), (4.2) with

$$\overline{v}(t,x) = f(w[t])(x) - f(u[t])(x) - \sum_{i=1}^m l_i(x) \varepsilon_i(t) + \tilde{v}(t,x) - v(t,x), \text{ for } (t,x) \in \Re_+ \times [0,1] \tag{4.30}$$

Therefore (4.4), (4.30), (2.7), (4.22) and (4.23) imply the following estimate for $t \geq 0$:

$$\|e\|_{[0,t]} \leq \sqrt{\max\left(|P|, \frac{Q}{2}\right)} \|e[0]\| + R \gamma \|e\|_{[0,t]} + \gamma \sum_{i=1}^m \|l_i\| \|\varepsilon_i\|_{[0,t]} + \gamma \|v - \tilde{v}\|_{[0,t]} \tag{4.31}$$

where $H(Q) := 2\sigma - \lambda_{N+1} - \sqrt{(2\sigma - \lambda_{N+1})^2 + 16Q^{-1}|L'PL|K^2}$, $\gamma := \sqrt{\dfrac{\tilde{g}}{2(\mu - \kappa)}}$, $\tilde{g} := \max\left(\dfrac{4|P|}{4\sigma + H(Q)}, \dfrac{Q}{2\lambda_{N+1}}\right)$ and

$\mu := (H(Q) + 2\lambda_{N+1})/4$. Combining estimates (4.29) and (4.31) and using (2.22), we obtain for all $t \geq 0$:

$$\|e\|_{[0,t]} \leq \sqrt{\max\left(|P|, \dfrac{Q}{2}\right)} \|e[0]\| + \gamma \exp(\kappa h) \sum_{i=1}^{m} \|l_i\| \|\xi_i\|_{[0,t]} + \Omega \|e\|_{[0,t]} + \left(1 + h\exp(\kappa h) \sum_{i=1}^{m} \|l_i\| \|c_i\|\right) \gamma \|v - \tilde{v}\|_{[0,t]} \quad (4.32)$$

Estimate (2.23) is a direct consequence of estimate (4.32) and (4.23). ◁

**Proof of Theorem 2.3:** Let (arbitrary) $u_0, w_0 \in D$, a bounded mapping $\xi : \Re_+ \to \Re^m$, a pair of inputs $v, \tilde{v} \in C^0(\Re_+ \times [0,1])$ that satisfy Assumption (H3) and an increasing sequence of times $\{t_j \geq 0 : j = 0,1,2,...\}$ with $t_0 = 0$, $\lim_{j \to +\infty}(t_j) = +\infty$ and $\sup_{j \geq 0}(t_{j+1} - t_j) \leq h$, where $h > 0$ satisfies (2.24). Consider the unique solution $u[t], w[t]$, $\zeta_i(t)$ ($i = 1,...,m$) of problem (2.4), (2.5), (2.10), (2.20), (2.16) with initial condition given by (2.9) and $w[0] = w_0$. Define for $t \geq 0$, $e[t]$ by (4.22) and

$$\varepsilon_i(t) := \int_0^1 c_i(x) e(t,x) dx \quad (i = 1,...,m). \quad (4.33)$$

Moreover, define $\|e\|_{[a,b]}, \|\xi_i\|_{[a,b]}, \|v - \tilde{v}\|_{[a,b]}$ by (4.23) for all $0 \leq a \leq b$. Let an arbitrary integer $j \geq 0$ be given. Using (2.17), (2.10), (2.20) and definitions (4.33), (4.22) we get for all $t \in (t_j, t_{j+1})$, $i = 1,...,m$:

$$\dot{\varepsilon}_i(t) = \int_0^1 \left(p \dfrac{d^2 c_i}{dx^2}(x) - q(x) c_i(x)\right) e(t,x) dx + \int_0^1 c_i(x)(f(w[t]))(x) dx - \int_0^1 c_i(x)(f(u[t]))(x) dx$$
$$+ \int_0^1 c_i(x)(\tilde{v}(t,x) - v(t,x)) dx - \sum_{r=1}^{m}\left(\int_0^1 c_i(x) l_r(x) dx\right) \xi_r(t_j) + \sum_{r=1}^{m}\left(\int_0^1 c_i(x) l_r(x) dx\right)\left(\int_0^1 k_r(s) e(t_j,s) ds\right) \quad (4.34)$$

It follows from (4.34), (2.7) and the Cauchy-Schwarz inequality that the following estimate holds for all $t \in (t_j, t_{j+1})$, $\tau \in (t_j, t]$:

$$|\dot{\varepsilon}(\tau)| \leq \left(\left\|p \dfrac{d^2 c_i}{dx^2} - q c_i\right\| + R\|c_i\|\right) \|e[\tau]\| + \|c_i\| \|\tilde{v}[\tau] - v[\tau]\| + \sum_{r=1}^{m}\left|\int_0^1 c_i(x) l_r(x) dx\right| |\xi_r(t_j)| + \|e[t_j]\| \sum_{r=1}^{m}\left|\int_0^1 c_i(x) l_r(x) dx\right| \|k_r\| \quad (4.35)$$

Using the facts that $t \in [t_j, t_{j+1})$, $\sup_{j \geq 0}(t_{j+1} - t_j) \leq h$ and (4.35), (4.23) we obtain the estimate:

$$|\varepsilon_i(t) - \varepsilon_i(t_j)| \exp(\kappa t) \leq h \exp(\kappa h) \left(\|c_i\| \|\tilde{v} - v\|_{[0,t]} + \sum_{r=1}^{m}\left|\int_0^1 c_i(x) l_r(x) dx\right| \|\xi_r\|_{[0,t]}\right)$$
$$\left(\left\|p \dfrac{d^2 c_i}{dx^2} - q c_i\right\| + R\|c_i\| + \sum_{r=1}^{m}\left|\int_0^1 c_i(x) l_r(x) dx\right| \|k_r\|\right) h \exp(\kappa h) \|e\|_{[0,t]} \quad (4.36)$$

Define for each $t \geq 0$:

$$\eta(t) := \max\{t_j : j = 0,1,2,..., t_j \leq t\} \quad (4.37)$$

Combining (4.36), (4.37) we obtain the following estimate for $t \geq 0$:

$$\sup_{0 \leq s \leq t}\left(|\varepsilon_i(s) - \varepsilon_i(\eta(s))| \exp(\kappa s)\right) \leq h \exp(\kappa h)\left(\|c_i\| \|\tilde{v} - v\|_{[0,t]} + \sum_{r=1}^{m}\left|\int_0^1 c_i(x) l_r(x) dx\right| \|\xi_r\|_{[0,t]}\right)$$
$$+ \left(\left\|p \dfrac{d^2 c_i}{dx^2} - q c_i\right\| + R\|c_i\| + \sum_{r=1}^{m}\left|\int_0^1 c_i(x) l_r(x) dx\right| \|k_r\|\right) h \exp(\kappa h) \|e\|_{[0,t]} \quad (4.38)$$

It follows from (2.4), (2.5), (2.6), (2.10), (2.16), (4.33), (4.22), (4.37) that $e[t]$ is a solution of (4.1), (4.2) with

$$\bar{v}(t,x) := f(w[t])(x) - f(u[t])(x) + \tilde{v}(t,x) - v(t,x)$$
$$+ \sum_{i=1}^{m} l_i(x) \left( \int_0^1 (k_i(s) - c_i(s)) e(\eta(t), s) ds + \varepsilon_i(\eta(t)) - \varepsilon_i(t) - \xi_i(\eta(t)) \right)$$
$$\text{for } (t,x) \in \Re_+ \times [0,1] \quad (4.39)$$

Therefore estimate (4.4) in conjunction with (4.39), (2.7), (4.37), (4.23) and the fact that $\sup_{j \geq 0}(t_{j+1} - t_j) \leq h$ imply the following estimate for $t \geq 0$:

$$\|e\|_{[0,t]} \leq \sqrt{\max\left(|P|, \frac{Q}{2}\right)} \|e[0]\| + \gamma \left( \exp(\kappa h) \sum_{i=1}^{m} \|l_i\| \|k_i - c_i\| + R \right) \|e\|_{[0,t]}$$
$$+ \gamma \left( \sum_{i=1}^{m} \|l_i\| \sup_{0 \leq s \leq t} (|\varepsilon_i(s) - \varepsilon_i(\eta(s))| \exp(\kappa s)) + \exp(\kappa h) \sum_{i=1}^{m} \|l_i\| \|\xi_i\|_{[0,t]} + \|v - \tilde{v}\|_{[0,t]} \right) \quad (4.40)$$

where $H(Q) := 2\sigma - \lambda_{N+1} - \sqrt{(2\sigma - \lambda_{N+1})^2 + 16 Q^{-1} |L'PL| K^2}$, $\gamma := \sqrt{\frac{\tilde{g}}{2(\mu - \kappa)}}$, $\tilde{g} := \max\left( \frac{4|P|}{4\sigma + H(Q)}, \frac{Q}{2\lambda_{N+1}} \right)$ and $\mu := (H(Q) + 2\lambda_{N+1})/4$. Combining estimates (4.38) and (4.30) and using (2.24), we obtain for all $t \geq 0$:

$$\|e\|_{[0,t]} \leq \sqrt{\max\left(|P|, \frac{Q}{2}\right)} \|e[0]\| + \left(1 + h \exp(\kappa h) \sum_{i=1}^{m} \|l_i\| \|c_i\| \right) \gamma \|v - \tilde{v}\|_{[0,t]}$$
$$+ \gamma \exp(\kappa h) \sum_{i=1}^{m} \left( \|l_i\| + \sum_{r=1}^{m} h \left| \int_0^1 c_r(x) l_i(x) dx \right| \|l_r\| \right) \|\xi_i\|_{[0,t]} + \Omega \|e\|_{[0,t]} \quad (4.41)$$

Estimate (2.25) is a direct consequence of estimate (4.41) and (4.23). ◁

## 5. Concluding Remarks

The present work showed that the extension of the approach in [17] to a wide class of 1-D, parabolic PDEs with non-local outputs is indeed feasible. Two different sampled-data observer designs are presented and analyzed: with and without an inter-sample predictor. Explicit conditions on the upper diameter of the (uncertain) sampling schedule were derived for both designs for exponential convergence of the observer error to zero in the absence of measurement noise and modeling errors. Moreover, explicit estimates of the convergence rate were deduced based on the knowledge of the upper diameter of the sampling schedule. On the other hand, when measurement noise and/or modeling errors are present, IOS estimates of the observer error were shown for both designs with respect to noise and modeling errors. Examples showed how the proposed methodology can be extended to other cases (e.g., boundary point measurements) and how we can detect cases of parabolic PDEs that allow exponential convergence of the observer error for arbitrarily large diameter of the sampling schedule for the sampled-data observer with the inter-sample predictor.

The use of the proposed sampled-data observers to sampled-data feedback stabilizers (see [8,11,18,28]) is straightforward. This is going to be the topic of future research. Another research direction that needs further study is the case of system of ODEs with hyperbolic PDEs as well as the case of 1-D, hyperbolic PDEs with non-local boundary conditions.